\documentclass[11pt]{article}
\usepackage[a4paper,margin=30mm]{geometry}
\usepackage{amsmath,amssymb,amsthm,mathtools}
\usepackage{enumitem}
\usepackage{microtype}
\usepackage[hidelinks]{hyperref}
\usepackage[T1]{fontenc}
\usepackage{lmodern}
\usepackage[utf8]{inputenc}

\newtheorem{theorem}{Theorem}[section]
\newtheorem{proposition}[theorem]{Proposition}

\newtheorem{corollary}[theorem]{Corollary}
\newtheorem{remark}[theorem]{Remark}
\newtheorem{example}[theorem]{Example}
\newtheorem{definition}[theorem]{Definition}

\newcommand{\R}{\mathbb{R}}
\newcommand{\C}{\mathbb{C}}
\newcommand{\norm}[1]{\left\lVert #1\right\rVert}
\newcommand{\ip}[2]{\left\langle #1,#2\right\rangle}

\title{Global Invertibility with Singular Points via the Mountain-Pass Principle}
\author{Dorota Bors\\
\small Faculty of Mathematics and Computer Science, University of Lodz\\[1ex]
Robert Sta\'nczy\\
\small Instytut Matematyczny, University of Wroc{\l}aw}
\date{October 2026}

\begin{document}
\maketitle

\begin{abstract}
Let $F:X\to H$ be a $C^1$ local homeomorphism from a real Banach space to a real Hilbert space. We prove that $F$ is a global homeomorphism if the functional $\Phi_y(x)=\tfrac12\norm{F(x)-y}^2$ satisfies the Palais--Smale condition for every $y\in H$. The proof uses Ekeland's variational principle and Hofer's characterization of mountain-pass points. At points where the first two derivatives vanish, a separation estimate for the cubic term yields a sufficient condition for local inversion. The complex cubic illustrates the role of this estimate. Applications include a nonlocal cubic equation on $L^2(0,1)$.
\end{abstract}

\section{Introduction}

Let $X$ be a real Banach space and $H$ a real Hilbert space. For a $C^1$ map $F:X\to H$ and $y\in H$, consider the functional
\begin{equation}\label{eq:phi}
\Phi_y(x)=\frac12\norm{F(x)-y}^2.
\end{equation}
The real map $x\mapsto x^3$ illustrates how an inverse can remain continuous at a singular value. For the complex cubic, the same formula gives three distinct preimages of each nonzero value. This distinction motivates a local separation condition on the cubic term.

\section{Branching of the complex cubic}

\begin{proposition}\label{prop:z3}
The map
\[
F:\R^2\simeq\C\longrightarrow\R^2\simeq\C,
\qquad F(z)=z^3,
\]
has an invertible derivative away from the origin. At the origin its first two derivatives vanish and its diagonal third derivative maps onto $\C$ with $\norm{F'''(0)[h,h,h]}=6\norm h^3$. Every $\Phi_y$ satisfies the Palais--Smale condition, and each nonzero value of $F$ has three preimages.
\end{proposition}

\begin{proof}
In real coordinates $z=x+iy$,
\[
F(x,y)=\bigl(x^3-3xy^2,\;3x^2y-y^3\bigr).
\]
Its derivative has determinant
\[
\det F'(x,y)=9(x^2+y^2)^2,
\]
so $F'(z)$ is invertible for $z\neq0$.  At the origin,
\[
F'(0)=F''(0)=0,
\qquad
F'''(0)[h,h,h]=6h^3.
\]
Hence for every $h$,
\[
\norm{F'''(0)[h,h,h]}=6\norm{h}^3,
\]
and the diagonal cubic $h\mapsto6h^3$ is onto $\C$.

For any $a\in\C$,
\[
\Phi_a(z)=\frac12|z^3-a|^2.
\]
If $\Phi_a(z_n)$ is bounded, then $|z_n|^3\le |z_n^3-a|+|a|$, so $(z_n)$ is bounded.  In finite dimension it has a convergent subsequence.  Thus $\Phi_a$ satisfies Palais--Smale.  Nevertheless
\[
1^3=\bigl(e^{2\pi i/3}\bigr)^3=\bigl(e^{4\pi i/3}\bigr)^3=1,
\]
so $F$ is three-to-one away from the origin in the usual angular sense and is not injective.
\end{proof}

\section{Global inversion under the Palais--Smale condition}

\begin{definition}
A $C^1$ functional $\varphi:X\to\R$ satisfies the Palais--Smale condition if every sequence $(x_n)$ for which $(\varphi(x_n))$ is bounded and $\norm{\varphi'(x_n)}\to0$ possesses a strongly convergent subsequence.
\end{definition}

We use the following geometric refinement of the mountain-pass principle due to Hofer \cite{Hofer}; see also Katriel \cite{Katriel}.

\begin{theorem}[Hofer's mountain-pass alternative, form used here]\label{thm:hofer}
Let $\psi\in C^1(X,\R)$ satisfy the Palais--Smale condition and the standard mountain-pass geometry between two points, and let $c$ be the associated minimax level.  Then the critical set at level $c$ contains either a local minimum or a mountain-pass point in Hofer's sense.  Thus one may choose $x_*\in X$ with $\psi(x_*)=c$ such that either $x_*$ is a local minimum or, for every neighborhood $N$ of $x_*$, the set
\[
\{x\in N:\psi(x)<\psi(x_*)\}
\]
is disconnected (and therefore nonempty).
\end{theorem}

The role of this strengthening is decisive: at a singular point the equation $\Phi_y'(x)=0$ no longer implies $F(x)=y$, but the \emph{local topology} of $\Phi_y$ still detects whether $F$ is locally one-to-one.

\begin{theorem}[Global inversion theorem]\label{thm:main}
Let $X$ be a real Banach space, $H$ a real Hilbert space, and let $F\in C^1(X,H)$.  Assume:
\begin{enumerate}[label=\textup{(\roman*)}]
\item $F$ is a local homeomorphism at every point of $X$;
\item for every $y\in H$, the functional
\[
\Phi_y(x)=\frac12\norm{F(x)-y}^2
\]
satisfies the Palais--Smale condition.
\end{enumerate}
Then $F:X\to H$ is a global homeomorphism.
\end{theorem}

\begin{proof}
We divide the proof into surjectivity and injectivity.

\medskip\noindent
\textbf{Surjectivity.}
Fix $y\in H$.  Since $\Phi_y\ge0$, Ekeland's variational principle yields a minimizing sequence $(x_n)$ such that
\[
\Phi_y(x_n)\longrightarrow m:=\inf_X\Phi_y,
\qquad
\norm{\Phi_y'(x_n)}\longrightarrow0.
\]
By Palais--Smale, after passage to a subsequence, $x_n\to x_*$.  Thus $x_*$ is a global minimizer of $\Phi_y$.

Suppose $F(x_*)\neq y$.  Since $F$ is a local homeomorphism at $x_*$, there are neighborhoods $U\ni x_*$ and $V\ni F(x_*)$ such that $F|_U:U\to V$ is a homeomorphism.  Put $z=F(x_*)$.  For all sufficiently small $t>0$,
\[
z_t=(1-t)z+ty\in V.
\]
Let $x_t=(F|_U)^{-1}(z_t)$.  Then
\[
\Phi_y(x_t)=\frac12\norm{z_t-y}^2
=(1-t)^2\frac12\norm{z-y}^2
<\Phi_y(x_*),
\]
contradicting minimality.  Therefore $F(x_*)=y$.  Since $y$ was arbitrary, $F$ is onto.

\medskip\noindent
\textbf{Injectivity.}
Suppose, to the contrary, that
\[
F(x_0)=F(x_1)=y,
\qquad x_0\neq x_1.
\]
We first establish genuine mountain-pass geometry without using a derivative estimate.  Choose a local chart $F|_U:U\to V$ at $x_0$.  Select $\rho>0$ so that
\[
\overline B_X(x_0,\rho)\subset U,
\qquad
\rho<\norm{x_1-x_0}.
\]
Let $G=(F|_U)^{-1}$.  By continuity of $G$ at $y$, there is $r>0$ such that
\[
B_H(y,r)\subset V,
\qquad
G(B_H(y,r))\subset B_X(x_0,\rho).
\]
If $\norm{x-x_0}=\rho$ and $\norm{F(x)-y}<r$, then $x=G(F(x))\in B_X(x_0,\rho)$, a contradiction.  Hence
\begin{equation}\label{eq:barrier}
\inf_{\norm{x-x_0}=\rho}\Phi_y(x)\ge \frac12 r^2>0.
\end{equation}
Since $\Phi_y(x_0)=\Phi_y(x_1)=0$, the translated functional has the standard mountain-pass geometry.  By the Palais--Smale condition and Theorem \ref{thm:hofer}, at its mountain-pass level $c\ge r^2/2$ there is a point $x_*$ which is either a local minimum or a mountain-pass point of $\Phi_y$.  In particular,
\[
\Phi_y(x_*)=c>0,
\qquad z:=F(x_*)\neq y.
\]

We show that neither alternative is possible.  Let $F|_{U_*}:U_*\to V_*$ be a local homeomorphism at $x_*$.  On $U_*$,
\[
\Phi_y=q_y\circ F,
\qquad
q_y(w)=\frac12\norm{w-y}^2.
\]
The point $z\neq y$ is not a local minimum of $q_y$: moving a short distance from $z$ toward $y$ strictly decreases $q_y$.  Hence $x_*$ is not a local minimum of $\Phi_y$.

Nor is $x_*$ a mountain-pass point.  Choose $\varepsilon>0$ so small that $B_H(z,\varepsilon)\subset V_*$.  With $R=\norm{z-y}$,
\[
\{w\in B_H(z,\varepsilon):q_y(w)<q_y(z)\}
=B_H(z,\varepsilon)\cap B_H(y,R).
\]
This set is nonempty and convex, hence path connected.  Its inverse image under the homeomorphism $F|_{U_*}$ is therefore path connected.  Thus, in the neighborhood $(F|_{U_*})^{-1}(B_H(z,\varepsilon))$, the strict lower sublevel of $\Phi_y$ is connected, contradicting the mountain-pass character of $x_*$.  This contradiction proves injectivity.

Finally, a bijective local homeomorphism is open, and hence its inverse is continuous.  Thus $F$ is a global homeomorphism.
\end{proof}

\begin{remark}\label{rem:regularcase}
If $F'(x)$ is an isomorphism for every $x$, the local-homeomorphism hypothesis follows from the Banach inverse mapping theorem.  Thus Theorem \ref{thm:main} contains the Hilbert-target Palais--Smale theorem of Idczak--Skowron--Walczak at the topological level.  Local homeomorphisms also allow singular derivatives.
\end{remark}

\section{A checkable cubic singularity criterion}

The abstract condition ``$F$ is a local homeomorphism'' can be verified at an odd-order singular point by higher-order inverse theorems.  We record the cubic case of Fija\l kowski's result \cite{Fijalkowski2003,Fijalkowski2018}.

\begin{theorem}[Fija\l kowski, cubic case]\label{thm:fij}
Let $X,Y$ be Banach spaces and $T\in C^3(U,Y)$ on a neighborhood $U$ of $x_0\in X$.  Suppose
\[
T'(x_0)=T''(x_0)=0,
\]
and define the homogeneous cubic
\[
P(h)=\frac16T'''(x_0)[h,h,h].
\]
Assume that $P:X\to Y$ is onto and that there exists $C>0$ such that
\begin{equation}\label{eq:fij}
C\norm{u}^2\norm{v}
\le
\norm{P(u+v)-P(u)}
\qquad (u,v\in X).
\end{equation}
Then $T$ is a local homeomorphism at $x_0$.
\end{theorem}

Condition \eqref{eq:fij} is a genuine nonbranching condition.  Taking $u=-h$ and $v=h$ gives
\[
\norm{P(h)}\ge C\norm{h}^3,
\]
so it implies the uniform cubic lower bound, but it is strictly stronger.

\begin{corollary}[Global inversion at cubic singularities]\label{cor:cubic}
Let $X$ be a real Banach space, $H$ a real Hilbert space, and $F\in C^3(X,H)$.  Assume that at each $x\in X$ one of the following holds:
\begin{enumerate}[label=\textup{(\alph*)}]
\item $F'(x):X\to H$ is a bounded linear isomorphism;
\item $F'(x)=F''(x)=0$ and, with
\[
P_x(h)=\frac16F'''(x)[h,h,h],
\]
$P_x$ is onto and there exists $C_x>0$ such that
\[
C_x\norm{u}^2\norm{v}
\le
\norm{P_x(u+v)-P_x(u)}
\quad\text{for all }u,v\in X.
\]
\end{enumerate}
If, in addition, $\Phi_y$ from \eqref{eq:phi} satisfies Palais--Smale for every $y\in H$, then $F$ is a global homeomorphism.
\end{corollary}

\begin{proof}
At points of type (a), $F$ is a local $C^1$-diffeomorphism by the inverse mapping theorem.  At points of type (b), Theorem \ref{thm:fij} gives a local homeomorphism.  Thus $F$ is a local homeomorphism everywhere, and Theorem \ref{thm:main} applies.
\end{proof}

\begin{remark}[Why $z^3$ is excluded]
For $P(z)=z^3$ on $\C$, choose $u=1$ and $u+v=e^{2\pi i/3}$.  Then $v\neq0$ but
\[
P(u+v)-P(u)=0.
\]
Thus the two distinct preimages violate the separation estimate \eqref{eq:fij}.
\end{remark}

\section{Examples}

\subsection{A scalar cubic singularity}

\begin{example}
Let
\[
F:\R\to\R,
\qquad F(x)=x^3+x^5.
\]
Then $F'(x)=3x^2+5x^4>0$ for $x\neq0$, while $F'(0)=F''(0)=0$ and
\[
P(h)=\frac16F'''(0)h^3=h^3.
\]
For $u,v\in\R$,
\[
|(u+v)^3-u^3|
=|v|\,|3u^2+3uv+v^2|
\ge \frac34 |u|^2|v|.
\]
Thus the cubic local inversion condition holds.  Moreover $|F(x)|\to\infty$ as $|x|\to\infty$, so in finite dimension every Palais--Smale sequence for $\Phi_y$ is bounded and has a convergent subsequence.  Corollary \ref{cor:cubic} therefore yields global invertibility.  Of course, here it can also be seen directly from strict monotonicity.
\end{example}

\subsection{A Hilbert-space cubic model with a singular derivative}

\begin{example}\label{ex:radial}
Let $H$ be any real Hilbert space and define
\begin{equation}\label{eq:radial}
F(x)=\bigl(\norm{x}^2+\norm{x}^4\bigr)x.
\end{equation}
Then $F$ is a global homeomorphism, $F'(0)=F''(0)=0$, and $F'(x)$ is an isomorphism for $x\neq0$.

The cubic part is
\[
P(x)=\norm{x}^2x.
\]
It is onto.  Furthermore
\[
DP(z)h=\norm{z}^2h+2\ip{z}{h}z.
\]
For $u,v\in H$,
\begin{align*}
\ip{P(u+v)-P(u)}{v}
&=\int_0^1\Bigl(\norm{u+tv}^2\norm{v}^2
+2\ip{u+tv}{v}^2\Bigr)\,dt\\
&\ge \norm{v}^2\int_0^1\norm{u+tv}^2\,dt.
\end{align*}
Since
\[
\int_0^1\norm{u+tv}^2dt
=\norm{u}^2+\ip{u}{v}+\frac13\norm{v}^2
\ge \frac14\norm{u}^2,
\]
we obtain
\[
\norm{P(u+v)-P(u)}\ge\frac14\norm{u}^2\norm{v}.
\]
Thus the higher-order local inversion criterion holds at zero.

For completeness, $\Phi_y$ satisfies Palais--Smale.  If $(\Phi_y(x_n))$ is bounded then $\norm{F(x_n)}$ is bounded, hence $(x_n)$ is bounded because
\[
\norm{F(x)}=\norm{x}^3+\norm{x}^5.
\]
After taking a subsequence, let $\norm{x_n}\to r$.  If $r=0$, then $x_n\to0$.  If $r>0$, the least spectral value of the positive self-adjoint operator $F'(x_n)$ is $\norm{x_n}^2+\norm{x_n}^4$, which is bounded below away from zero.  Consequently $\norm{\Phi_y'(x_n)}\to0$ implies $F(x_n)-y\to0$.  Since the radial inverse of \eqref{eq:radial} is continuous, $x_n$ converges strongly.  Hence Theorem \ref{thm:main} applies even in infinite-dimensional $H$.
\end{example}

This example also illustrates why ``critical point of $\Phi_y$'' is the wrong substitute for ``solution of $F(x)=y$'' at a singular derivative: for every $y$, one has $\Phi_y'(0)=0$, although $F(0)=y$ only when $y=0$.

\subsection{A polynomial with a rank-one derivative}

Consider the polynomial map
\[
F(a,b)=\bigl(a^3+a^5-b^5,\;b^3+b^5+a^5\bigr).
\]
Its Jacobian is
\[
DF(a,b)=
\begin{pmatrix}
3a^2+5a^4 & -5b^4\\
5a^4 & 3b^2+5b^4
\end{pmatrix}.
\]
At $(a,0)$ with $a\ne0$, the derivative has rank one. Corollary~\ref{cor:cubic} therefore gives no conclusion about local inversion at these points.

\section{A nonlocal cubic equation}

Consider the integral operator
\begin{equation}\label{eq:origint}
F(x)=A(x^2)x+r(x),
\qquad
(Az)(t)=\int_0^1K(t,s)z(s)\,ds,
\end{equation}
where $K\in C([0,1]^2)$ and $0<\alpha\le K(t,s)\le\beta$. The operator $A$ maps bounded subsets of $L^1(0,1)$ into uniformly bounded, equicontinuous subsets of $C([0,1])$ and is compact by the Arzel\`a--Ascoli theorem. Its associated cubic term is
\[
P(h)=A(h^2)h.
\]
The pointwise positivity of $K$ yields
\[
\norm{P(h)}_{L^2}\ge \alpha\norm{h}_{L^2}^3,
\]
Proposition~\ref{prop:z3} illustrates why one also needs a condition on local injectivity. For the following class of kernels, the inverse can be obtained explicitly.

\subsection{A weighted cubic on $L^2(0,1)$}

\begin{proposition}\label{prop:weighted}
Let $H=L^2(0,1)$ and let $k\in L^\infty(0,1)$ satisfy
\[
0<\alpha\le k(s)\le\beta<\infty\quad\text{a.e.}
\]
Set
\[
q(x)=\int_0^1k(s)x(s)^2\,ds,
\qquad
F(x)(t)=q(x)x(t).
\]
Then $F:H\to H$ is a global homeomorphism.  It has a cubic singularity at $0$ and satisfies the Fija\l kowski increment condition
\[
\norm{F(u+v)-F(u)}_2
\ge\frac{\alpha^2}{4\beta}\norm{u}_2^2\norm{v}_2.
\]
For every $y\in H$, $\Phi_y$ satisfies Palais--Smale.
\end{proposition}

\begin{proof}
Define the equivalent Hilbert inner product
\[
B(u,v)=\int_0^1k(s)u(s)v(s)\,ds,
\qquad q(u)=B(u,u).
\]
For $P=F$,
\[
DP(z)h=q(z)h+2B(z,h)z.
\]
Therefore
\begin{align*}
B(P(u+v)-P(u),v)
&=\int_0^1\bigl(q(u+tv)q(v)+2B(u+tv,v)^2\bigr)\,dt\\
&\ge q(v)\int_0^1 q(u+tv)\,dt\\
&\ge \frac14 q(u)q(v).
\end{align*}
The last inequality is the same one-dimensional quadratic estimate used in Example \ref{ex:radial}, but in the $B$-norm.  Since
\[
|B(w,v)|\le\beta\norm{w}_2\norm{v}_2,
\qquad
q(w)\ge\alpha\norm{w}_2^2,
\]
the claimed increment bound follows.

Surjectivity is explicit.  If $y\neq0$, put
\[
x=q(y)^{-1/3}y.
\]
Then $q(x)x=y$; for $y=0$ take $x=0$.  The increment estimate implies injectivity and local invertibility, and the explicit formula gives continuity of the inverse.

To verify Palais--Smale, note first that $\norm{F(x)}_2=q(x)\norm{x}_2\ge\alpha\norm{x}_2^3$, so a bounded $\Phi_y$-sequence is norm bounded.  Along a subsequence let $\norm{x_n}_2\to r$.  If $r=0$, there is strong convergence to zero.  If $r>0$, $DF(x_n)$ is an isomorphism and its inverse norm is bounded uniformly on that subsequence (this follows either from the rank-one formula for $DF$ or from the positive lower bound in the $B$-inner product).  Since
\[
\Phi_y'(x_n)=DF(x_n)^*(F(x_n)-y)\to0,
\]
we obtain $F(x_n)\to y$, and the continuous inverse gives strong convergence of $x_n$.
\end{proof}

This is the case of \eqref{eq:origint} with $r=0$ and a kernel independent of $t$, namely $K(t,s)=k(s)$.

\section{An implicit-function consequence}

Applying Theorem~\ref{thm:main} at each fixed parameter gives the following solvability result.

\begin{corollary}[Global solvability for parameter slices]\label{cor:implicit}
Let $P$ be a parameter set and let $G:P\times X\to H$.  Suppose that for each fixed $p\in P$, the slice
\[
G_p:X\to H,\qquad G_p(x)=G(p,x),
\]
is $C^1$, is a local homeomorphism, and for every $y\in H$ the functional
\[
x\longmapsto\frac12\norm{G(p,x)-y}^2
\]
satisfies Palais--Smale.  Then for every $(p,y)\in P\times H$ there exists a unique $x=S(p,y)$ such that
\[
G(p,x)=y.
\]
For each fixed $p$, the map $y\mapsto S(p,y)$ is continuous.  At points where $D_xG(p,x)$ is an isomorphism, the usual implicit-function theorem gives the corresponding $C^1$ local dependence on parameters.  At a genuinely cubic singularity one should in general expect only topological, not $C^1$, inversion; the scalar cube root is the basic model.
\end{corollary}

Joint continuity in $p$ requires an additional uniform compactness/properness hypothesis on the family; it does not follow merely from slice-wise invertibility.

\section{Relation with the literature}

Katriel developed mountain-pass criteria for global homeomorphisms in topological and metric settings \cite{Katriel}.  Idczak, Skowron and Walczak proved the smooth Hilbert-target result for local $C^1$-diffeomorphisms under Palais--Smale conditions for the squared-distance functionals \cite{IdczakSkowronWalczak}.  Gut\'u later generalized that line to Banach targets using the Chang Palais--Smale condition for local $C^1$-diffeomorphisms \cite{Gutu}.  Fija\l kowski established higher-order local inversion at odd singularities \cite{Fijalkowski2003} and a different global inversion theorem with singular points based on Hadamard--L\'evy type control \cite{Fijalkowski2018}.

The formulation of Theorem \ref{thm:main} is a natural synthesis of these ingredients: it replaces ``local $C^1$-diffeomorphism'' by ``$C^1$ local homeomorphism'', thereby allowing singular derivatives, while retaining the ordinary Hilbert-space Palais--Smale functional.  Notably, Gut\'u explicitly observed that Katriel's proof suggests precisely a conjectural principle of the form ``local homeomorphism + a Palais--Smale condition on distance to $y$ implies bijectivity'' \cite{Gutu}.  The proof above realizes that principle using Hofer's topological description of the mountain-pass critical point.

\section*{Acknowledgments}
We thank Darek Idczak (University of Lodz) for first drawing our attention to the injectivity issue. We used ChatGPT to identify the complex cubic counterexample, correct the mathematical argument, and typeset the manuscript in \LaTeX.

\end{document}